\documentclass{amsart}

\theoremstyle{plain}
\newtheorem{theorem}{Theorem}[section]

\newcommand{\PGL}{\operatorname{PGL}}

\newcommand{\gxg}{G\times G}
\newcommand{\bxb}{B\times B^-}

\newcommand{\diag}{\operatorname{diag}}

\newcommand{\End}{\operatorname{End}}

\newcommand{\bA}{{\mathbb A}}
\newcommand{\bC}{{\mathbb C}}

\newcommand{\bP}{{\mathbb P}}

\newcommand{\bR}{{\mathbb R}}

\newcommand{\bZ}{{\mathbb Z}}

\newcommand{\cF}{{\mathcal F}}

\newcommand{\cL}{{\mathcal L}}
\newcommand{\cM}{{\mathcal M}}

\newcommand{\cO}{{\mathcal O}}

\newcommand{\cT}{{\mathcal T}}

\newcommand{\cV}{{\mathcal V}}

\newcommand{\inv}{^{-1}}

\newcommand{\wP}{\widetilde P}

\DeclareMathSymbol{\curvearrowright}{\mathrel}{AMSb}{"79}
\DeclareMathSymbol\rightsquigarrow {\mathrel}{AMSa}{"20}
\DeclareMathSymbol\square {\mathord}{AMSa}{"03}
\DeclareMathSymbol{\ltimes}         {\mathbin}{AMSb}{"6E}
\DeclareMathSymbol{\nmid}           {\mathrel}{AMSb}{"2D}
\DeclareMathSymbol{\twoheadrightarrow}  {\mathrel}{AMSa}{"10}

\newcommand{\ratmap}{- \kern -3pt \to}

\newcommand{\Cone}{\operatorname{Cone}}

\newcommand{\const}{\operatorname{const}}
\newcommand{\Vol}{\operatorname{Vol}}

\newcommand{\Proj}{\operatorname{Proj}}

\newcommand{\Aut}{\operatorname{Aut}}

\theoremstyle{plain}

\newtheorem{lemma}[theorem]{Lemma}

\newtheorem{proposition}[theorem]{Proposition}
\newtheorem*{theorem-no}{Theorem}

\newtheorem{conjecture}[theorem]{Conjecture}

\theoremstyle{definition}

\newtheorem{definition}[theorem]{Definition}

\newtheorem{notations}[theorem]{Notations}         

\newtheorem{example}[theorem]{Example}

\newtheorem{remark}[theorem]{Remark}

\newtheorem*{acknowledgements}{Acknowledgments} 

\theoremstyle{remark}

\date{December 27, 2003}
\title{On K-Stability of Reductive Varieties}
\author{Valery Alexeev and Ludmil Katzarkov}
\address{Department of Mathematics\\
University of Georgia\\
Athens, GA 30602, USA}
\email{valery@math.uga.edu}
\address{Department of Mathematics\\
UC Irvine \\
Irvine, CA, 92612, USA}
\email{lkatzark@math.uci.edu}

\newcommand{\Dvert}{D_{\rm vert}}
\newcommand{\Dhor}{D_{\rm hor}}
\newcommand{\kxk}{K\times K}
\newcommand{\LambdaR}{\Lambda_{\bR}}
\newcommand{\dbd}{\bar{\partial}\partial}

\begin{document}
\bibliographystyle{amsalpha}
\maketitle

\begin{abstract}
  G. Tian and S.K. Donaldson formulated a conjecture relating
  GIT stability of a polarized algebraic variety to the existence of a
  K\"ahler metric of constant scalar curvature.  In
  \cite{Donaldson_Stability} Donaldson partially confirmed it
  in the case of projective toric varieties. In this paper
  we extend Donaldson's results and computations to a new case, that
  of reductive varieties.
\end{abstract}

\section*{Introduction}

Around 1997 G. Tian and S.K. Donaldson formulated a conjecture
relating GIT stability of a polarized algebraic variety to the
existence of a K\"ahler metric of constant scalar curvature, see
\cite{Tian_KE,BC,SS} and
\cite{Donaldson_FieldsLecture,AS,Donaldson_Stability}.  The general
idea of such a relationship has been known as a ``folklore
conjecture'' for a long time. It comes naturally from earlier
works of Yau and Tian.

Here, we will refer to the Donaldson's version of the conjecture
(\ref{conj:Donaldson} below), since we will make a strong use of
notations and results developed in \cite{Donaldson_Stability}. In this
beautiful paper Donaldson partially confirmed conjecture
\ref{conj:Donaldson} in the case of projective toric varieties, where
his arguments are based on the following foundation blocks:
\begin{enumerate}
\item a very well-known correspondence between projective toric
  varieties and lattice polytopes (see
  e.g.  \cite{Fulton_ToricVarieties,Oda_ConvexBodies});
\item a somewhat lesser-known correspondence between convex piecewise
  linear functions on polytopes and equivariant degenerations of
  projective toric varieties (see e.g. 
  \cite{KapranovSturmfelsZelevinsky_ChowPolytopes92} or
  \cite{Alexeev_CMAV}); and
\item the general framework developed by Guillemin and Abreu
  \cite{Guillemin_KaehlerStructures,Abreu_ExtremeMetrics} which
  relates invariant K\"ahler metrics on projective toric varieties to
  convex functions on the corresponding polytopes.
\end{enumerate}

Now let $G$ be a complex reductive group with Weyl group $W$. Then to
every $W$-invariant maximal-dimensional lattice polytope $P$ in the
weight lattice of a maximal subtorus of $G$ one can, in a rather
elementary way, associate an equivariant projective normal
compactification $V_{P}$ of $G$, generalizing the correspondence~(1)
above.

On the other hand, in
\cite{AlexeevBrion_Affine,AlexeevBrion_Projective} the first author
and M. Brion built a theory of degenerations parallel to (2), and
obtained an answer that is formally very similar: Any $W$-invariant
rational convex PL function on $P$ defines a $W$-invariant subdivision
of $P$ and a degenerating family $\cV\to \bC$ in which every fiber
$\cV_t$ with $t\ne0$ is isomorphic to $V$ and the special fiber
$\cV_0$ is a \emph{stable reductive variety} in the sense of
\cite{AlexeevBrion_Affine,AlexeevBrion_Projective}.

The purpose of this paper is to extend the Guillemin-Abreu-Donaldson
theory and results and computations of \cite{Donaldson_Stability} to
this new and much larger class of projective $G$-compactifications and
reductive varieties.  In particular, our Theorem~\ref{thm:main} is a
generalization of one of the main results of
\cite{Donaldson_Stability}. We use Theorem 3.7 together with  Donaldson's
examples \cite{Donaldson_Stability} of polarized toric surfaces that do not
admit metrics of constant scalar curvature to produce new
\emph{non-toric} examples. Instead of 2-dimensional toric, our examples  are
8-dimensional reductive varieties.

The paper is organized as follows: In Section 2 we review the  relevant
definitions and results from \cite{Donaldson_Stability}.  In Section 3
we recall some results about reductive varieties from
\cite{AlexeevBrion_Affine,AlexeevBrion_Projective}. In Section 4 we
extend several of Donaldson's results from \cite{Donaldson_Stability}
%(e.g. Proposition 4.2.1, Proposition 3.2.4, Lemma 3.2.6 ) 
to the case of reductive varieties.  Section 5 contains examples of
reductive varieties which do not admit K\"ahler metrics of constant
scalar curvature.

\begin{acknowledgements}
  We are very grateful to M. Brion, S.K. Donaldson and W.~Graham for
  many helpful discussions, and to G. Tian and the referee for useful
  suggestions. The second author also would like to thank V. Apostolov
  for explaining Abreu's work. S. Donaldson informed us that his student
  A. Raza has studied some differential geometric aspects of the
  constant scalar curvature problem for reductive varieties. 
 
  It is our pleasure to thank IHES, where a large part of this work
  was done for wonderful working conditions.  The first author's work
  was partially supported by NSF grant 0101280.
\end{acknowledgements}

Throughout the paper, we will use the following

\begin{notations}
  $K$ will denote a compact real reductive group of dimension $n$ with
  a maximal subtorus $T$ of dimension $r$. We will denote by
  $G=K_{\bC}$ and $H=T_{\bC}$ their complexifications.  We choose a
  Borel subgroup $B$ of $G$ and an opposite Borel subgroup $B^-$ so
  that $B\cap B^-=H$.
  
  Let $\Lambda\simeq\bZ^r$ denote the group of characters of $H$. It
  comes with an action of the Weyl group $W$ and with a decomposition
  of $\LambdaR$ into Weyl chambers. $\LambdaR^+$ will denote the
  positive chamber, and $r_i\in \Lambda$ the simple roots.
\end{notations}

\section{Basic definitions and results from \cite{Donaldson_Stability}} 

The precise formulation of the conjecture relating K-stability and
existence of K\"ahler metric of constant scalar curvature was given by
Donaldson in \cite{Donaldson_Stability}:

\begin{conjecture}\label{conj:Donaldson} 
  A smooth polarized projective variety $(V,L)$ admits a K\"ahler metric
  of constant scalar curvature in $c_{1}(L)$ if and only if it is
  K-stable.
\end{conjecture}
 
The definition of K-stability \cite[Def.2.1.2]{Donaldson_Stability}
involves another space $V_{0}$ which is allowed to be a general
scheme. Let $\cF$ be an ample line bundle over a projective scheme
$W$, and suppose one has a fixed $\bC^{*}$-action on the pair
$(W,\cF)$. For each positive integer $k$ one has a vector space
$$
H_{k} = H^{0}(W, \cF^{k})
$$
with a $\bC^{*}$-action. From this, one obtains integers $d_{k}=
\text{dim}\ H_{k}$ and $w_{k}$, the weight of the induced action on
the highest exterior power of $H_{k}$. The integers $d_{k},w_{k}$ are,
for large $k$, given by polynomial functions of $k$, with rational
coefficients: $d_{k}=Q(k), w_{k}=P(k)$ say. Define $F(k)= w_{k}/k
d_{k}$.  For large enough $k$ one has an expansion
$$
F(k)= F_{0} + F_{1} k^{-1} + F_{2} k^{-2} + \dots, $$
with rational
coefficients $F_{i}$. The \emph{Futaki invariant} of the
$\bC^{*}$-action on $(W,\cF)$ is defined to be the coefficient
$F_{1}$.
 
\begin{definition}
  A \emph{test configuration} or \emph{test family} for $(V,L)$ of
  exponent $r$ consists of
  \begin{enumerate}
  \item a scheme ${\mathcal V}$ with a $\bC^{*}$-action;
  \item a $\bC^{*}$-equivariant line bundle ${\mathcal L}\rightarrow
    {\mathcal V}$;
  \item a flat $\bC^{*}$-equivariant map $\pi:{\mathcal V}\rightarrow
    \bC$, where $\bC^{*}$ acts on $\bC$ by multiplication in the
    standard way;
\end{enumerate}
such that any fiber $V_{t}= \pi^{-1}(t)$ for $t\neq 0$ is isomorphic
to $V$ and the pair $(V,L^{r})$ is isomorphic to $(V_{t}, {\mathcal
  L}\vert_{V_{t}})$.
\end{definition} 
             
\begin{definition}
  The pair $(V,L)$ is \emph{K-stable} if for each test configuration
  for $(V,L)$ the Futaki invariant of the induced action on
  $(V_{0},{\mathcal L}\vert_{V_{0}})$ is less than or equal to zero,
  with equality if and only if the configuration is a product
  configuration.
  
  The pair $(V,L)$ is \emph{equivariantly K-stable} if one restricts
  oneself only to equivariant test configurations. For a toric variety
  this means $H$-invariant families, and in our case this will
  mean $\gxg$-equivariant families.
\end{definition} 

See also the discussion in \cite{Donaldson_Stability} showing that
K-stability is closely related to assymptotic GIT stability.
             
Following \cite{Donaldson_Stability} we briefly recall some facts
about Mabuchi functional.  It is a real-valued function $\cM$ on the
set of K\"ahler metrics in the same K\"ahler class $[\omega_{0}]$,
defined up to the addition of an overall constant. The metrics of
constant scalar curvature are critical for the Mabuchi functional. The
functional is defined through the formula for its variation at a
metric $\omega= \omega_{0} +i\dbd \psi$ with respect to an
infinitesimal change $\delta \psi$ in the K\"ahler potential:

$$   
\delta {\mathcal M} = \int_{V} (S-a) \delta \psi \frac{\omega^{n}}{n!}
$$ 
  
Here, $S$ is the scalar curvature of $\omega$ and $a$ is the average
value of the scalar curvature. Therefore $\delta {\mathcal M}$ is not
changed if one adds a constant to $\delta \psi$, and so depends only
on the variation of the metric. 

\begin{definition}
  A functional $I_V$ on the set of K\"ahler potentials is defined by the
  formula
  $$ \delta I_V = 2\int_V \delta\psi \frac{\omega^n}{n!} $$
\end{definition}

Assume that $D=\sum D_i$ is an effective divisor on $V$ representing
$c_1(V)$, with smooth components $D_i$ (such $D$ exists for any smooth
toric variety $V$). Let $\chi$ be a meromorphic form with $(\chi) =
-D$ and let $\nu = |\chi|^{-2}$.

\begin{proposition}[\cite{Donaldson_Stability}, 3.2.4]
  \label{prop:linear_part}
  For any metric $\omega= \omega_0 + 2i\dbd\psi$ on $V$
  $$
  \cM(\omega) = \cL_a(\omega) + \int_V \log\nu \frac{\omega^n}{n!},
  \text{ where } 
  \cL_a(\omega) = -I_D(\psi) + a I_V(\psi)$$ is the linear part
  of Mabuchi functional.
\end{proposition}

Now, let $(V,L)$ be a polarized projective toric variety with a moment
polytope $P$. Via the Guillemin-Abreu theory, an equivariant K\"ahler
potential on $V$ corresponds to a convex function $u$ on the interior
of $P$ with a logarithmic asymptotic behavior near the
boundary. Namely, if faces are defined by linear inequalities
$\delta_k (x)\ge 0$ then
$$
u = \sum \delta_k(x)\log \delta_k(x) + 
\text{ a smooth function.}
$$

The lattice determines standard measures $d\mu$ on the
polytope $P$ and $d\sigma$ on faces of~$P$.
In \cite[3.2.6]{Donaldson_Stability} Donaldson establishes the
following formula for the linear part of Mabuchi functional:
$$ 
\cL_a =  (2\pi)^{n} \left(
\int_{\partial P} u\, d\sigma - a \int_P u\, d\mu 
\right)
$$
On the other hand, in \cite[4.2.1]{Donaldson_Stability} Donaldson
shows that Futaki invariant of an equivariant test configuration
defined by a rational PL function $f$ is given by a very similar
formula
$$ 
-F_1 = \frac{1}{2\Vol P} \left(
\int_{\partial P} f\, d\sigma - a \int_P f\, d\mu 
\right)
$$
Thus, one obtains an ingenious connection between two seemingly
very different invariants. Donaldson uses this
connection to prove, among other things, the following
\begin{proposition} [\cite{Donaldson_Stability}, Prop. 7.1.2]
  Suppose there is a function $f\in C^1$ with $\cL_a(f) <0 $. Then the
  Mabuchi functional is not bounded below on the invariant metrics and
  the manifold $V$ does not admit any K\"ahler metric of constant
  scalar curvature in the given cohomology class.
\end{proposition}

Our aim will be to generalize these results to the reductive
case.

%%%%%%%%%%%%%%%%%%%%%%%%%%%%%%%%%%%%%%%%%
\section{Overview of reductive varieties}
\label{sec:Overview of reductive varieties}

\subsection{Complexes of polytopes and (stable) reductive varieties}

We recall some of the results and constructions of
\cite{AlexeevBrion_Projective}.

Lattice points $\lambda\in\Lambda^+$ are in bijection with 
irreducible $G$-representation $E_{\lambda}$. The algebra of regular
functions on $G$ can be written canonically as 
\begin{eqnarray*}
&&\bC[G] = \oplus_{\lambda\in\Lambda} \End E_{\lambda},
\text{ and }
 \End E_{\lambda} \cdot \End E_{\mu} \subset
\oplus \End E_{\lambda+\mu - \sum_{n_i\ge0} n_ir_i}.
\end{eqnarray*}

Let $P\subset \Lambda$ be a maximal-dimensional $W$-invariant
polytope with vertices in $\Lambda$.  Let $P^+ = P\cap
\LambdaR^+$ be the part of $P$ lying in the positive chamber, and
$\Cone\Delta^+ \subset \bR\oplus\LambdaR$ be the cone over
$(1,P^+)$. The vector space
$$ R_{P} = \oplus_{\lambda\in\bZ\oplus\Lambda\cap\Cone\Delta^+} 
\End E_{\lambda} $$
has a natural structure of a subalgebra in $\bC[\bC^*\times G]$. It is
finitely generated e.g. by \cite[4.8]{AlexeevBrion_Affine}.

One defines $V_P = \Proj R_P$ and $L_P= \cO(1)$. Then $V_P$ is a
normal projective equivariant $\gxg$-compactification of group $G$
(with $\gxg$ acting on $G$ by left and right multiplication:
$(g_1,g_2).g = g_1\inv g g_2$), and $L_P$ is a $\gxg$-linearized ample
sheaf on $V_P$.

The fixed point set $(V^{\diag H},L^{\diag H})$ is a toric variety
with a $H$--linearized ample line bundle. It corresponds to the same
polytope $P$, and it comes with an action by $WH$, a semidirect
product of $W$ and $H$. The $\gxg$-orbits of $V$ are in bijection
with $WH$-orbits of $V^{\diag H}$ and with $W$-orbits of faces of
the polytope $P$.

\begin{example}
 Take $G=\PGL_n$ embedded into $V=\bP Mat_{n\times n} =
\bP^{n^2-1}$. In this case, $\LambdaR$ is $\bR^{n-1}$ divided into
$n!$ Weyl chambers. The polytope $P$ is a simplex with a vertex on one
of the rays of the positive chamber $\LambdaR^+$, and the other $n-1$
vertices are its reflections under $W=S_n$.
\end{example}

\begin{example}
In the previous example, take $n=2$. The polytope $P$ is an interval
symmetric about the origin.  The variety $V=\bP^3$ is the wonderful
compactification of De Concini-Procesi for $\PGL_2$. 
\end{example}

\begin{example}
  For any semisimple group $G$, consider a point in the interior of
  $\LambdaR^+$ and let $P$ be the convex hull of its
  $W$-reflections. The corresponding variety $V_P$ is called a
  wonderful compactification of $G$. 

  For $G=\PGL_3$ and $W=S_3$, $P$ is a hexagon and $P^+$ is a
  4-gon. Note that the  wonderful compactification of $\PGL_3$ is not
  toric. 
\end{example}

\cite[Sec.5]{AlexeevBrion_Affine} and
\cite[Sec.2]{AlexeevBrion_Projective} generalize this picture to the
case of a $W$-invariant complex of polytopes, i.e. a finite
$W$-invariant collection $\Delta=\{P_i\}$ of lattice polytopes in
$\LambdaR$ such that the intersection $P_i\cap P_j$ of any two
polytopes is a union of faces of both. Note that the polytopes need
not be maximal-dimensional and that reductive varieties in general do
not contain $G$.

\begin{theorem}[\cite{AlexeevBrion_Projective}, Thm. 2.8]
  \label{thm:srv-classification}
  \begin{enumerate}
  \item A complex of polytopes $\Delta$ defines a family
    of polarized stable reductive varieties $\{ (V_{\Delta,t},
    L_{\Delta,t}) \}$ parameterized by a certain cohomology group
    $H^1(\Delta/W, \Aut)$. The choice $t=1$ gives a distinguished
    ``untwisted'' member of this family.
  \item $V^{\diag H}_{\Delta,t}$ is a stable toric variety, a union of
    ordinary toric varieties corresponding to polytopes $P_i$ in
    $\Delta$, together with a $W$-action.
  \item The $\gxg$-orbits of $V_{\Delta,t}$ are in a bijection with
    $WH$-orbits of $V^{\diag H}_{\Delta,t}$ and with the set
    $\Delta/W$.
  \item Variety $V_{\Delta,t}$ is irreducible iff $\Delta/W$ contains
    a unique maximal polytope. In this case, $H^1(\Delta/W,
    \Aut)=\{1\}$. Variety $V_{\Delta}$ is a $G$-compactification iff
    $\Delta$ consists of one $W$-invariant polytope $P$ of maximal
    dimension plus its faces.
  \end{enumerate}  
\end{theorem}

Note that here we use the notion of stable reductive varieties as
defined in \cite{AlexeevBrion_Projective} and not the notion of
K-stability. 

\subsection{Nonsingular reductive varieties}

In toric geometry, a lattice polytope $P$ corresponds to a projective
toric variety $V_P$ together with a linearized ample sheaf $L_P$. It
is well-known that the  variety $V_P$ is non-singular if and only if $P$ is
a \emph{Delzant polytope}, i.e. at every vertex precisely $\dim P$
edges meet and the integral generators of these edges form a basis of
the lattice. We will need the following generalization of this result
to the case of reductive varieties, which we will use in the case of a
single $W$-invariant polytope.

\begin{proposition}\label{prop:nonsingularity}
  Let $\Delta=\{P_i\}$ be a $W$-invariant complex of polytopes, and
  $V=V_{\Delta}$ be the corresponding (stable) reductive variety. 
  \begin{enumerate}
  \item If $V$ is nonsingular, then $\Delta$ is a disjoint union of
    Delzant polytopes plus its faces.
  \item If $\Delta$ is a disjoint union of Delzant polytopes none of
    whose vertices lie on the supporting hyperplanes of $\LambdaR^+$,
    plus its faces, then $V$ is nonsingular.
  \end{enumerate}
\end{proposition}
\begin{proof}
  If $H$ is a torus acting on a smooth variety $Z$, then the fixed
  point set $Z^H$ is also nonsingular. By
  Theorem~\ref{thm:srv-classification}, $V^{\diag H}$ is a stable toric
  variety, a union of ordinary toric varieties corresponding to
  $P_i\in\Delta$. It is nonsingular if and only if the maximal
  polytopes are disjoint and each of them is Delzant.

  In the opposite direction, if none of the polytopes $P_i$ contain
  vertices on the supporting hyperplanes of the positive chamber, then the 
  variety $V$ is toroidal, i.e. locally analytically it is isomorphic
  to $\bA^N$ times the toric variety corresponding to the polytopes
  $P_i$. Hence, it is nonsingular if $P_i$ are Delzant.  
\end{proof}

%%%%%%%%%%%%%%%
\subsection{Canonical class}

According to \cite{Brion_CurvesDivs}, the anticanonical divisor of any
spherical variety for group $G$ can be written as
$$
-K_V = \partial_G V + \partial_B G
$$
Here $\partial_G V=\sum D_v$ is the reduced sum of $G$-invariant
divisors, and $\partial_B G=\sum n_{\rho}D_{\rho}$ is the union of
$B$-invariant, non-$G$-invariant divisors (``colors''), with uniquely
defined positive coefficients $n_{\rho}$. For spherical varieties in
general the coefficients $n_{\rho}$ can be arbitrary, and
\cite{Brion_CurvesDivs} provides the precise formula. However, for
reductive varieties, which are spherical for the  group $\gxg$, all
coefficients $n_{\rho} =2$, see
\cite[Sec.5.2]{AlexeevBrion_Projective}. We will call $\partial_{\gxg}
V$ the \emph{vertical} and $\partial_{\bxb} V$ the \emph{horizontal}
parts of the anticanonical divisor (this notation comes from the
behavior of these boundaries under the moment map).  The above
formula then becomes
$$ 
-K_V = \Dvert + \Dhor
$$ 

If $\Delta$ is a complex of polytopes satisfying the conditions of
Proposition \ref{prop:nonsingularity}(2) then all codimension one
faces also satisfy \ref{prop:nonsingularity}(2). Hence, every $D_v$ is
a smooth reductive variety. In addition, since $V$ is toroidal, the
linear system $|\Dhor|$ is basepoint-free. Hence, there exists a
smooth divisor representing $\Dhor$.

\subsection{One-parameter degenerations}

Let $f$ be a convex rational $W$--invariant PL function on $P$. 
In the same fashion as in the toric case, 
$f$ defines a $(\dim V +1)$-dimensional ``test'' family $\cV\to \bC$
such that every fiber $\cV_t$ for $t\ne 0$ is isomorphic to $V$. The
construction, contained in \cite[Sec. 4.2]{AlexeevBrion_Projective},
is as follows. 

After replacing the polytope $P$ by a large divisible multiple $NP$,
i.e. replacing $L$ by $L^N$, one can assume that the domains of
linearity of $f$ are lattice polytopes $P_i$.

Consider a bigger polytope $\wP$ in $\bR\oplus \LambdaR$ which is
bounded from below by $(0,P)$ and from above by the graph of $R-f$ for
some $R\gg0$.  Then, by the above construction, the variety $V_{\wP}$
is a projective compactification of the group $\bC^* \times G$, and it
comes with a natural map $\pi$ to $\bP^1$. The test family $\cV$ is
$V_{\wP}\setminus \pi\inv{(\infty)}$. The special fiber $\pi\inv(0)$
is a stable reductive variety for the complex of polytopes
$\Delta=\{P_i\}$.

%%%%%%%%%%%%%%%%%%%%%%%%%%%%%%%%%%%%%%%%%%%%%%%%%%%%%%%%%%%%%%%%%%%%%%
\section{Linear part of Mabuchi functional and Futaki invariant}

Every symplectic variety $V$ with a Hamiltonian action of a compact
Lie group $K$ admits a moment map $m:V\to k^*$ to the dual of the Lie
algebra. The moment map commutes with the $K$--action on $V$ and with
the coadjoint action of $K$ on $ k^*$.  Therefore, the image is a
union of coadjoint orbits. 

Let $t \subset k$ be the Cartan subalgebra and let $t^* \subset k^*$
be a splitting of the projection $k^*\to t^*$. Since every coadjoint
orbit intersects the positive Weyl chamber $(t^*)^+$ at one point, one
obtains a continuous map $\pi: V\to (t^*)^+=\LambdaR^+$.  By a theorem
of Guillemin-Sternberg \cite{GuilleminSternberg_MM2} the image of $V$
under the moment map $m_T$ for a maximal torus $T\subset K$ is a
polytope, and by a result of Kirwan \cite{Kirwan_MM3} the image
$\pi(V)$ is a polytope. Both polytopes live in the space $\LambdaR$
and both are traditionally called moment polytopes. For us, $\pi(V)$
will be important.

In the case when $V$ is a projective spherical variety for
the complexified Lie group $G=K_{\bC}$, some general facts about the
moment map were established in \cite{Brion_MM} via representation
theory. In particular, $\pi(V)$ identifies with the Brion's moment
polytope of a spherical variety.  For a reductive variety
corresponding to a polytope $P$ one thus obtains the moment polytope
$P^+$.

Let $\lambda\in \LambdaR^+$ be an integral vector.  It corresponds to
an irreducible $G$--representation $E_{\lambda}$. According to the  Weyl
character formula, $\dim E_{\lambda}$ is given by a polynomial
$$
h(x) = h_e(\lambda) + h_{e-1}(\lambda) + \dots 
$$
written as a sum of its homogeneous parts.

Any reductive variety is spherical for the action of $\gxg$. For each
$\lambda\in \Lambda^+$ the corresponding $\gxg$-representation is
$\End E_{\lambda}$.  Let us write its dimension as the sum of its
homogeneous parts
$$ 
\dim \left(\End E_{\lambda}\right) =  H(\lambda) = 
h^2(\lambda) = H_d(\lambda) + H_{d-1}(\lambda) + \dots
$$ 
and extend each $H_d$ to a polynomial function on $\LambdaR^+$.

In the previous section, we wrote $-K_V$ and the $B\times
B^-$-boundary of $V$ in the form $\Dvert + \Dhor$. The vertical part
$\Dvert$ is a union of $\gxg$-invariant divisors, each of them is a
reductive variety corresponding to a face of polytope $P$ modulo the
$W$-action. Under the map $\pi$, $\Dvert$ maps to a union of faces of
$P^+$. As in the toric case, for every face $F\subset P$ we will
denote by $d\sigma$ the Euclidean measure on $F$ induced by the
sublattice $\Lambda_F = \Lambda\cap \bR F$.

The horizontal part $\Dhor$ is not $\gxg$-invariant, and it is easy to
see that under the map $\pi$ it maps surjectively onto $P^+$.

\begin{lemma}\label{volume2}
  \begin{enumerate}
  \item The push-forward of the Liouville measure on $V$ is 
    $$
    \pi_*\frac{\omega^n}{n!} = (2\pi)^r H_d(x) \,d\mu
    $$
  \item Similarly, the push-forward of the Liouville measure on
    $\Dvert$ is 
    $$
    \pi_*\frac{\omega^{n-1}}{(n-1)!} = 
    (2\pi)^r H_d(x) \,d\sigma.
    $$
  \item The push-forward of the Liouville measure on 
    $\Dhor$ is 
    $$\pi_*\frac{\omega^{n-1}}{(n-1)!} = (2\pi)^r\cdot 2H_{d-1}(x) 
    \,d\mu.$$
  \end{enumerate}
\end{lemma}
\begin{proof}
  This comes out straight from the Riemann-Roch theorem applied to $H^0(
  L^s)$ on the variety $V$ and on a single coadjoint orbit
  $\cO_{\lambda}$ for $s\gg 0$. The push-forward of the
  Liouville measure is given by integrating $\omega^n/n!$. On the
  other hand, the Riemann-Roch theorem gives
  $$ 
  H^0(L^s) = \int ch(L^s) \cdot \operatorname{Td}(\cT) =
  s^n \frac{c_1(L)^n}{n!} 
  -s^{n-1} \frac{1}{2} K \frac{c_1(L)^{n-1}}{(n-1)!}  + \dots
  $$
  So, the volume of a single coadjoint orbit $\cO_{\lambda}$ is
  given by $H_d(\lambda)$ and the volume of the restriction of the
  horizontal part of $-K_V$ to $m\inv(\cO_{\lambda})$ by
  $2H_{d-1}(\lambda)$.  Reductive varieties are spherical for
  $\gxg$-action, so in the decomposition of $H^0(V,L^s)$ every $\End
  E_{\lambda}$, $\lambda\in sP^+$ appears with multiplicity one; and
  the formula follows.
\end{proof}

\begin{proposition}
  $\kxk$--invariant K\"ahler metrics on $V$ are in a bijection with
  $W$--invariant symplectic convex potentials $u$ on $P$ which have
  the same behavior near the boundary of $P$ as in the toric case
  (i.e. $u=\sum l_i \log l_i + $ a smooth function).
\end{proposition}
\begin{proof}
  By \cite[Thm.3.1]{GuilleminSternberg_MM2}, $\omega= i
  \partial\bar\partial \phi $ for a some potential $\phi$ on $V$. The
  restriction $\phi_T$ of this potential to the toric variety
  $V^{\diag T}$ is $W$--invariant, and hence by the toric case,
  applying the Legendre transform, gives a symplectic potential $u$ on
  $P^+$. Clearly, it has to be $W$--invariant.
  
  Vice versa, every $W$--invariant symplectic potential $u$ on $P$
  gives a $W$--invariant potential $\phi_T$ on $V^{\diag T}$. By
  $\kxk$--action, it extends to a unique potential $\phi$ on $V$.
  Here, we are using the fact that for a spherical variety the
  preimages of coadjoint orbits under the moment map are precisely the
  $\kxk$--orbits.
\end{proof}

\begin{theorem}\label{thm:Futaki}
  Let $f$ be a convex rational $W$--invariant PL function on $P$. Then
  the Futaki invariant of the corresponding test family is given by
  the formula
  $$ -F_1(f)= 
  \frac{1}{2 \int_{P^+} H_d}
  \left( 
    \int_{\partial P^+} f H_d \,d\sigma + 
    2 \int_{P^+} f H_{d-1}\, d\mu 
    - a \int_{P^+} f H_d \,d\mu 
  \right),
  $$
  where
  $$ a= \frac{\int_{\partial P^+} H_d \,d\sigma +
    2\int_{P^+} H_{d-1} 
    \,d\mu}{\int_{P^+} H_d \,d\mu}
  $$
\end{theorem}
\begin{proof}
  Same as the proof of Proposition 4.2.1 of
  \cite{Donaldson_Stability}. The difference is that this time an
  integral point $\lambda\in P^+$ represents not a one-dimensional
  vector space, as in the toric case, but the vector space $\End
  E_{\lambda}$ whose dimension is $H(\lambda)$. We use the next lemma
  to estimate the sum, and get
\begin{eqnarray*}
  A&=& \int_{P^+} (R-f) H_d d\mu \\
  B&=& \frac12 \int_{\partial^+} (R-f) H_d d\sigma
        + \int_{P^+} (R-f) H_{d-1} d\mu \\ 
  C&=& \int_{P^+} H_d d\mu \\
  D&=& \frac12 \int_{\partial P^+} H_d d\sigma
        + \int_{P^+} H_{d-1} d\mu 
\end{eqnarray*}
Substituting $F_1=(AD-BC)/2C^2$ gives the formula.
\end{proof}

\begin{remark}
  In the toric case one has $H_d=1$, $H_{d-1}=0$, $P^+ = P$ and the
  formula reduces to that of \cite[4.2.1]{Donaldson_Stability}.
\end{remark}

\begin{lemma}
  Let $P\subset \bR^n$ be a lattice polytope of dimension $n$ and
  $H_d:\bR^n\to \bR$ be a homogeneous polynomial function of degree
  $d$. Then for $k\gg0$,
$$ \sum_{\lambda\in kP\cap\bZ^n} H_d(\lambda) =
        k^{n+d} \int_P H_d d\mu +
        \frac12 k^{n+d-1} \int_{\partial P} H_d d\sigma
        + O(k^{n+d-2})
$$
\end{lemma}
\begin{proof}
  Elementary. For a monomial $x_1^{a_1}\dots x_n^{a_n}$ the
  computation easily reduces to the case of a polytope of dimension
  $n+d$ and the function $H_0=1$, i.e. to counting integral points in
  a polytope, where the analogous formula is well-known
  (cf. \cite[A1]{Donaldson_Stability}).
  Alternatively, this is a corollary of a general formula of Pukhlikov
  and Khovanski \cite{Pukhlikov_Khovanski}.
\end{proof}

\begin{theorem}\label{thm:linpartMabuchi}
  Let $u$ be a $W$--invariant symplectic potential on $P$
  corresponding to a K\"ahler form $\omega$. Then the linear part of
  the Mabuchi functional is given by the formula
  $$ {\mathcal L}_a(u) = 
  (2\pi)^r
  \left( 
    \int_{\partial P^+} uH_d \,d\sigma + 
    2 \int_{P^+} uH_{d-1} \,d\mu 
   -a  \int_{P^+} u H_d \,d\mu
  \right),
  $$
  where
  $$ a= \frac{\int_{\partial P^+}H_d \,d\sigma +
    2\int_{P^+} H_{d-1} 
    \,d\mu}{\int_{P^+} H_d \,d\mu}    
  $$
\end{theorem}
\begin{proof}
  As in the proof of Proposition 3.2.4 and Lemma 3.2.6 of
  \cite{Donaldson_Stability}, we use the formula from Proposition
  \ref{prop:linear_part} and compute
  \begin{eqnarray*}
  I_D(\psi) - a I_V (\psi) 
  &=& I_{\Dvert}(\psi) + I_{\Dhor}(\psi) - a I_V (\psi) \\
  &=& \int_{\Dvert}\psi\frac{\omega^{n-1}}{(n-1)!} +
  \int_{\Dhor}\psi\frac{\omega^{n-1}}{(n-1)!} - 
  a \int_{\Dhor}\psi\frac{\omega^n}{n!} 
  \end{eqnarray*}
  Applying Lemma \ref{volume2} completes the computation.
\end{proof}

Therefore, for equivariant compactifications of $G$, similar to the
toric case, the two seemingly unrelated functions, Futaki invariant of
an equivariant test configuration and the linear part of Mabuchi
functional, are given by the same formula. Although they are defined
on different sets, convex PL functions in the first case and convex
$C^{\infty}$--functions with a prescribed boundary behavior in the
second, we can approximate one class of functions by another one
freely. Hence, in the cases when the sign of $F_1$, resp. ${\mathcal
  L}_a$, is definite, we get the equivalence between (equivariant)
K--stability and boundedness of Mabuchi functional from below. 
As a corollary, we obtain

\begin{theorem} [cf. \cite{Donaldson_Stability}, Prop. 7.1.2]
  \label{thm:main}
  Let $V$ be a smooth equivariant compactification of a complex
  reductive group $G$.  Suppose there is a function $f\in C^1$ with
  $\cL_a(f) <0 $. Then the Mabuchi functional is not bounded below on
  the invariant metrics and the manifold $V$ does not admit any
  K\"ahler metric of constant scalar curvature in the given cohomology
  class.
\end{theorem}

\begin{remark}
  This theorem also applies, with the same proof, to a larger class of
  smooth reductive varieties described by
  Proposition~\ref{prop:nonsingularity}(2). 
\end{remark}

\section{Examples of $G$-compactifications without CSC
 metrics}

We modify Example \cite[7.2]{Donaldson_Stability} appropriately to get
a new series of varieties $V_n$ without a K\"ahler metric of constant
scalar curvature. 

Donaldson starts with a triangle $(0,0)$, $(0,1)$,
$(1,0)$. 
Working with the corner $A=(0,0)$ first, he considers the part of the
first quadrant bounded by the $x$- and $y$-axes and points $B=
(0,1/4)$, $C=(1/4,0)$ and $D_n =(r_n,r_n)$, where $(n-2)/4(3n-5)$. He
then smoothes out the corners $B$ and $C$ by adding very short
segments of various slopes so that the resulting polyhedral body is
Delzant, i.e. at every vertex the integral generators of edges form a
basis of the standard lattice $\bZ^2$.

He then translates this picture symmetrically to the other three
corners of the original triangle, and calls the resulting polytope
$P_n$. It is clear that $P_n$ is a Delzant polytope which is a
``smoothing'' of a 9-gon, and that the corresponding smooth toric
variety $V_n$ can be obtained by repeating blow-ups of the projective
plane.  The polytope represents an ample sheaf $L_n$ on $V_n$.
Donaldson shows that $P_n$ does not contain a metric of constant
curvature in the cohomology class $c_1(L_n)$.

We take $G=PGL_3$. Then Weyl group is $S_3$ acting on $\Lambda=\bZ^2$
and $\Lambda_{\bR}$ is divided into 6 cones, the chambers.  Take a
point $A$ far in the interior of the positive chamber, at the same
distance from both sides. Reflect it symmetrically to obtain vertices
of a hexagon. At each of the 6 vertices, repeat Donaldson's
construction close to the corners.

The resulting polytope $P_n$ is then a Delzant polytope
approximating a polytope with $3\cdot 6 = 18$ vertices.  The
corresponding polarized variety $V_n$ is a compactification of $G$ and
has dimension 8.  It is nonsingular by
Proposition~\ref{prop:nonsingularity}.

The Weyl character formula gives:
$$ H_d = x^2y^2(x+y)^2, \qquad 
   H_{d-1}= xy(x+y)^3 + 2 x^2 y^2 (x+y) $$
Here, $(x,y)$ are natural coordinates on $\LambdaR^+\cap \Lambda \simeq
\bZ_{\ge0}^2$. 

We claim that $P_n$ does not have a metric of constant curvature in
the cohomology class $c_1(L_n)$ by exhibiting a convex function $f$ on
it with $\cL_a(f)<0$. We use the formula of
Theorem~\ref{thm:linpartMabuchi} to compute it.

First, we need to observe that the function $2H_{d-1}/H_d$ reaches its
minimum near the corner $A$, therefore near the corner $D_n$ of $P^+$
the function $2H_d-aH_d$ is negative.

Define the function $f=f_{\epsilon,n}\ge0$ to be a convex function which
is equal to zero on the inside of $P_n$ and has a simple crease at the
distance $\epsilon$ to the point $D_n$. Define it at other corners by
symmetry, to be $W$-invariant. We claim that for $n\gg 1/\epsilon \gg
0$ one has ${\mathcal L}_a(f)<0$, as required. Indeed, the expressions

$$ 
\int_{P^+} H_d \,d\mu, \qquad 
\int_{\partial P^+} H_d \,d\sigma, \qquad
\int_{P^+} H_{d-1}\,d\mu 
$$
all have finite positive limits near $D_n$. The contribution of
$\partial P^+$ to $\cL_a(f)$ in Theorem~\ref{thm:linpartMabuchi} is
$O(\epsilon/n)$, as Donaldson's computation shows, and is positive.
Finally, the contribution of the integral over $P^+$ is approximately
$
\epsilon^2 (2H_{d-1}-aH_d)\,d\mu.
$
As we observed, this is negative 
and grows as $\const\cdot\epsilon^2$. Since $\epsilon^2\gg
\epsilon/n$, we get $\cL_a(f)<0$.
We note that the  varieties $V_n$ are not toric.  They are, however,
reductive varieties.

In conclusion, we note that Donaldson \cite{Donaldson_Stability}
includes a number of conjectures and tentative results in the toric
case. Using Theorems \ref{thm:Futaki} and \ref{thm:linpartMabuchi} we
can generalize them to our situation. Once they are established in the
toric case, we expect that they can be proved in the reductive case,
as well.

%\bibliography{string,primary11-15-03}
\providecommand{\bysame}{\leavevmode\hbox to3em{\hrulefill}\thinspace}
\providecommand{\MR}{\relax\ifhmode\unskip\space\fi MR }
\providecommand{\MRhref}[2]{%
  \href{http://www.ams.org/mathscinet-getitem?mr=#1}{#2}
}
\providecommand{\href}[2]{#2}

\end{document}